\documentclass[a4paper,12pt]{amsart}
\usepackage{amsmath}
\usepackage{amsfonts}
\usepackage{amssymb}
\usepackage{hyperref}
\usepackage{fourier}
\usepackage[a4paper,left=2.5cm,right=2.5cm,top=2.5cm, bottom=2cm]{geometry}
\usepackage[ddmmyyyy]{datetime}
\usepackage{tikz}
\usepackage{lineno}
\usepackage{float}
\usepackage{pagecolor}
\newtheorem{lemma}{Lemma}%[section]
\newtheorem{theorem}[lemma]{Theorem}
\newtheorem{corollary}[lemma]{Corollary}

\theoremstyle{definition}\newtheorem{example}[lemma]{Example}

\newenvironment{proofof}{\noindent}{\hfill$\Box$\medskip}

\newcommand{\cont}{\subseteq}

\newcommand{\s}{^*}\newcommand{\del}{\backslash}

\newcommand{\defin}{\textbf}

\title{A Note on Construction of Dual-Hamiltonian Graphs}
\author[Jo\~ao Paulo Costalonga]{Jo\~ao Paulo Costalonga$^1$}
\thanks{$^1$Departamento de Matem\'atica, Universidade Federal do Esp\'irito Santo. Av. Fernando Ferrari, 514; Campus de Goiabeiras,
  29075-910, Vit\'oria, ES, Brazil. e-mail: {\upshape joaocostalonga@gmail.com}.}
\begin{document}

\begin{abstract}
A connected simple graph is said dual-hamiltonian if its vertex set has a $2$-coloring such that each color class induces a tree. We call such a coloring a hamiltonian coloring. We prove that if $G$ is a graph with a certain type of hamiltonian coloring and $T$ is a tree, then $G\times T$ is also dual-hamiltonian having the same certain type of hamiltonian coloring. This result is used to constructed a class of dual-hamiltonian graphs, which includes the hypercubes and other multidimensional grids.
\end{abstract}

\maketitle

A \defin{bond} in a connected graph is a minimal set with the property that its deletion results in a disconnected graph. A bond in a connected graph $G$ has at most $\|G\|-|G|+2$ edges; when a bond has this size it is a \defin{hamiltonian bond}. A graph is \defin{dual-hamiltonian} if it has a hamiltonian bond. Some authors call dual-hamiltonian graphs \defin{Yutsis graphs}. An alternative definition for dual-hamiltonian graphs is provided by the following result:

\begin{theorem}(Jaeger~\cite{Jaeger})\label{criterion} A subset $B\cont E(G)$ is a hamiltonian bond of a connected graph $G$ if and only if $G\del B$ has two induced trees of $G$ as connected components.
\end{theorem}

If $G$ is a planar graph with dual graph $G\s$, then each hamiltonian bond of $G$ corresponds to a hamiltonian cycle of $G^*$. Determining if a graph is dual-hamiltonian is an NP-complete problem. Indeed, it is in the NP class by the certificate provided by Theorem\ref{criterion} and it is NP-complete because determining if a planar graph is hamiltonian is an NP-complete problem~\cite{Garey}.

The most notorious problem about dual-hamiltonian graphs is probably Jaeger's conjecture~\cite{Jaeger} that all cyclically $4$-connected cubic graphs are dual-hamiltonian. Whitney~\cite{Whitney} proved that the planar case of this conjecture is equivalent to the Four Colors Theorem. Another motivation to study dual-hamiltonian graphs is their application in quantum mechanics (see \cite{Van Dyck} for more references).

Some necessary conditions for dual-hamiltonianicity are established in \cite{Aldred}. We can add one more that follows from Theorem \ref{criterion} with a straightforward calculation:

\begin{corollary}
If $G$ is a dual-hamiltonian graph, then each induced subgraph $H$ of $G$ satisfies $\|H\|\le|H|-2+|H|^2/4$.
\end{corollary}

To the best of our knowledge, in opposition to hamiltonianicity, not many sufficient conditions for dual-hamiltonianicity are established in literature. This make establishing if some important classes of graphs are dual-hamiltonian a bit hard. Our objective here is to provide some tools for constructing dual-hamiltonian graphs using cartesian product of graphs. Our results were first motivated by the conjecture of Haidong Wu and Melissa Flynn (private communication) that all hypercubes are dual-hamiltonian; the $n$-dimensional hypercube $Q_n$ may be defined as the graph with $\{0,1\}^n$ as vertex-set and such that two vertices are adjacent if and only if they differ by exactly one coordinate. We will establish that Flynn's and Wu's conjecture is true in Corollary \ref{hypercube}.

Before stating the main result, some definitions are necessary. Let $G$ be a graph and let $C$ be a $2$-coloring of $V(G)$. We say that $C$ is a \defin{hamiltonian coloring} of $G$ if each color class of $C$ induces a tree in $G$. By Theorem \ref{criterion}, a $2$-coloring $C$ of $G$ is hamiltonian if and only if the set of edges with edvertices of different colors is a hamiltonian bond of $G$. If $X\cont V(G)$ and $C$ is a $2$-coloring of $V(G)$, we denote by $C\Delta X$ the $2$-coloring obtained from $C$ by switching the colors of the vertices in $X$. If $C$ is a hamiltonian coloring of $G$ and $I$ and $J$ are disjoint $2$-subsets of $V(G)$, we say that $(I,J)$ is a \defin{quartet} of $C$ if:
\begin{enumerate}
  \item [{\bf (Q1)}] both $I$ and $J$ have one vertex of each color,
  \item [{\bf (Q2)}] $C\Delta I$ is a hamiltonian coloring of $G$, and
  \item [{\bf (Q3)}] each color class of $C\Delta J$ induces a forest with exactly two connected components: one meeting $I$ and the other meeting $J$.
\end{enumerate}

The \defin{cartesian product} $G\times H$ of two graphs $G$ and $H$ is defined as the graph with vertex set $V(G)\times V(H)$ and edge-set as the union of $\{(u,x)(v,x):uv\in E(G)$ and $x\in V(H)\}$ and $\{(u,x)(u,y): u\in V(G)$ and $xy\in E(H)\}$. The next theorem is our main result:

\begin{theorem}\label{tree}
Let $G$ be a connected graph and $T$ be a tree. If $G$ has a hamiltonian coloring with a quartet, then $G\times T$ also has a hamiltonian coloring with a quartet.
\end{theorem}

The variation of Theorem \ref{tree} where the  hamiltonian colorings do not necessarily have quartets is not valid. Indeed, let $T_1$ and $T_2$ be disjoint trees, each one having a simple path with three edges, and let $G$ be the graph obtained from $T_1\cup T_2$ by adding all possible edges linking a vertex of $T_1$ to a vertex of $T_2$. Note that $\{V(T_1),V(T_2)\}$ is the unique $2$-partition of $V(G)$ such that both members induce forests in $G$. This implies that $G\times T$ is not dual-hamiltonian for each tree $T$ with at least two vertices. This paper contains only one more section, which is dedicated to proving Theorem \ref{tree}. 

From Theorem \ref{tree} we derive the next corollary. This result was proved for $n\le 5$ and conjectured in general in an unpublished undergraduate thesis written by Melissa Flynn under supervision of Haidong Wu.

\begin{corollary}\label{hypercube}
For $n\ge 2$, the $n$-dimensional hypercube $Q_n$ is dual-hamiltonian. Moreover, for $n\ge 3$, $Q_n$ has a hamiltonian coloring with a quartet.
\end{corollary}

We denote by $P_n$ the path graph with $n$ edges. It is straightforward to check that $Q_1$ and $Q_2$ are dual-hamiltonian. Since $Q_n\cong P_1\times\cdots\times P_1$, the second part of Corollary \ref{hypercube} follows by induction using Theorem \ref{tree} with the next example as initial case.

\begin{example}\label{ex-cube}
Let $n\ge 1$. The coloring of $P_1\times P_1\times P_n$ in Figure \ref{fig1} is hamiltonian with $\big((i_1,i_2),(j_1,j_2)\big)$ as quartet. In particular this proves that the $3$-dimensional cube has a hamiltonian coloring with a quartet.

\begin{figure}[h]\centering
\begin{tikzpicture}
\tikzstyle{gray} = [shape = circle,fill = gray,minimum size = 7pt,inner sep=1pt]
\tikzstyle{black} =[shape = circle,fill = black,minimum size = 7pt,inner sep=1pt]
%   C
%
%D
%
%
%    B
%A
\node[gray]  (a1) at (0,0) {};
\node[gray]  (b1) at (0.5,0.5) {};
\node[gray]  (c1) at (0.5,2.5) {};
\node[black] (d1) at (0,2) {};

\node[black] (a2) at (2,0) {};
\node[black] (b2) at (2.5,0.5) {};
\node[gray ] (c2) at (2.5,2.5) {};
\node[black] (d2) at (2,2) {};

\node[black] (a3) at (4,0) {};
\node[gray]  (b3) at (4.5,0.5) {};
\node[gray ] (c3) at (4.5,2.5) {};
\node[gray]  (d3) at (4,2) {};

\node[black] (a4) at (6,0) {};
\node[black] (b4) at (6.5,0.5) {};
\node[gray ] (c4) at (6.5,2.5) {};
\node[black] (d4) at (6,2) {};

\node[black] (a5) at (8,0) {};
\node[gray]  (b5) at (8.5,0.5) {};
\node[gray]  (c5) at (8.5,2.5) {};
\node[gray]  (d5) at (8,2) {};

\node[black] (a6) at (10,0) {};
\node[black] (b6) at (10.5,0.5) {};
\node[gray ] (c6) at (10.5,2.5) {};
\node[black] (d6) at (10,2) {};

\node() at (10.5,0) {$\cdots$};
\node() at (11,0.5) {$\cdots$};
\node() at (11,2.5) {$\cdots$};
\node() at (10.5,2) {$\cdots$};

\tikzstyle{gline} = [draw = gray!110, line width = 1.7pt]
\tikzstyle{bline} = [draw = black, line width = 1.7pt]
\tikzstyle{tline} = [densely dotted, line cap=round]

\draw (d1)--(d2)--(d3)--(d4)--(d5)--(d6);
\draw (b1)--(b2)--(b3)--(b4)--(b5)--(b6);
\draw (a1)--(a2);
\draw (a1)--(b1)--(c1)--(d1)--(a1);
\draw (a2)--(b2)--(c2)--(d2)--(a2);
\draw (a3)--(b3)--(c3)--(d3)--(a3);
\draw (a4)--(b4)--(c4)--(d4)--(a4);
\draw (a5)--(b5)--(c5)--(d5)--(a5);
\draw (a6)--(b6)--(c6)--(d6)--(a6);

\draw[gline](a1)--(b1)--(c1)--(c2)--(c3)--(c4)--(c5)--(c6);
\draw[gline](d3)--(c3)--(b3);
\draw[gline](d5)--(c5)--(b5);

\draw[bline](d1)--(d2)--(a2)--(a3)--(a4)--(a5)--(a6);
\draw[bline](a2)--(b2);
\draw[bline](d4)--(a4)--(b4);
\draw[bline](d6)--(a6)--(b6);

\node() at (0.5,2.8) {$j_1$}; %c1
\node() at (2.5,2.8) {$i_1$}; %c2
\node() at (-0.3,2) {$j_2$};%d1
\node() at (2.75,0.75) {$i_2$};%b2
\end{tikzpicture}
\caption{The figure for Example \ref{ex-cube}}\label{fig1}
\end{figure}
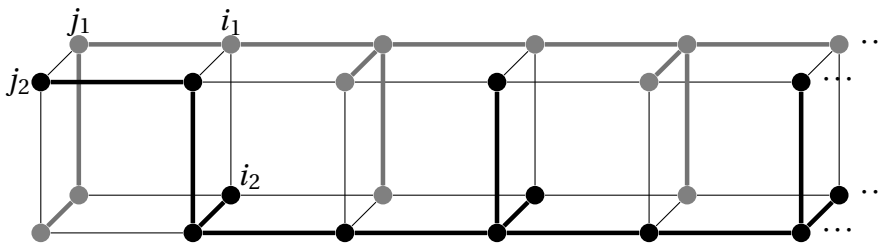
\end{example}

Similarly to the class of hypercubes, using the next example we may prove that a larger class of multidimensional grids are dual-hamiltonian using the graphs in th next example as initial cases. But we do not know if all the multidimensional grids are dual-hamiltonian. We conjecture so, but we conjecture that some do not have hamiltonian colorings with quartets.

\begin{example}
For $m\ge2$ and $n\ge 3$, the grid $P_m\times P_n$ has a hamiltonian coloring with a quartet $\big((i_1,i_2),(j_1,j_2)\big)$. The coloring is exemplified in Figure \ref{fig-even-grid} for odd values of $n$ and in Figure \ref{fig-odd-grid} for even values.

\begin{figure}[h]\centering
\begin{minipage}{6.5cm}\centering
\begin{tikzpicture}
\tikzstyle{gray}  = [shape = circle, fill = gray!110, draw = gray!110, very thick, minimum size = 2pt,inner sep=2pt]
\tikzstyle{black} = [shape = circle, fill = black, draw = black, very thick , minimum size = 2pt,inner sep=2pt]

\node[gray](11)at(0,5){};\node[black](12)at(1,5){};\node[black](13)at(2,5){};\node[black](14)at(3,5){};\node[black](15)at(4,5){};
\node[gray](21)at(0,4){};\node[gray] (22)at(1,4){};\node[gray] (23)at(2,4){};\node[gray] (24)at(3,4){};\node[black](25)at(4,4){};
\node[gray](31)at(0,3){};\node[black](32)at(1,3){};\node[black](33)at(2,3){};\node[black](34)at(3,3){};\node[black](35)at(4,3){};
\node[gray](41)at(0,2){};\node[gray] (42)at(1,2){};\node[gray] (43)at(2,2){};\node[gray] (44)at(3,2){};\node[black](45)at(4,2){};
\node[gray](51)at(0,1){};\node[black](52)at(1,1){};\node[black](53)at(2,1){};\node[black](54)at(3,1){};\node[black](55)at(4,1){};
\node[gray](61)at(0,0){};\node[gray] (62)at(1,0){};\node[gray] (63)at(2,0){};\node[gray] (64)at(3,0){};\node[black](65)at(4,0){};
\pgfdeclarelayer{bg}    % declare background layer
\pgfsetlayers{bg,main} 
\begin{pgfonlayer}{bg}
\draw (11)--(15);\draw (21)--(25);\draw (31)--(35);\draw (41)--(45);\draw (51)--(55);\draw (61)--(65);
\draw (11)--(61);\draw (12)--(62);\draw (13)--(63);\draw (14)--(64);\draw (15)--(65);
\end{pgfonlayer}
\node at (-0.3,5){$i_1$};
\node at (-0.3,3){$j_1$};
\node at (4.3,4) {$j_2$};
\node at (4.3,0) {$i_2$};

\tikzstyle{gline} = [draw = gray!110, line width = 1.7pt]
\tikzstyle{bline} = [draw = black, line width = 1.7pt]
\tikzstyle{tline} = [densely dotted, line cap=round]
\draw[bline](12)--(13)--(14)--(15)--(25)--(35)--(34)--(33)--(32);
\draw[bline](35)--(45)--(55)--(65);
\draw[bline](55)--(54)--(53)--(52);
\draw[gline](11)--(21)--(22)--(23)--(24);
\draw[gline](21)--(31)--(41)--(42)--(43)--(44);
\draw[gline](41)--(51)--(61)--(62)--(63)--(64);
\end{tikzpicture}
\caption{A hamiltonian coloring with a quartet for $P_5\times P_6$.}\label{fig-even-grid}
\end{minipage}
\begin{minipage}{6.5cm}\centering
\begin{tikzpicture}
\tikzstyle{gray}  = [shape = circle, fill = gray!110, draw = gray!110, very thick, minimum size = 2pt,inner sep=2pt]
\tikzstyle{black} = [shape = circle, fill = black, draw = black, very thick , minimum size = 2pt,inner sep=2pt]

\node[black](11)at(0,6){};\node[black](12)at(1,6){};\node[black](13)at(2,6){};\node[black](14)at(3,6){};\node[black](15)at(4,6){};
\node[gray] (21)at(0,5){};\node[gray] (22)at(1,5){};\node[gray] (23)at(2,5){};\node[gray] (24)at(3,5){};\node[black](25)at(4,5){};
\node[gray] (31)at(0,4){};\node[black](32)at(1,4){};\node[black](33)at(2,4){};\node[black](34)at(3,4){};\node[black](35)at(4,4){};
\node[gray] (41)at(0,3){};\node[gray] (42)at(1,3){};\node[gray] (43)at(2,3){};\node[gray] (44)at(3,3){};\node[black](45)at(4,3){};
\node[gray] (51)at(0,2){};\node[black](52)at(1,2){};\node[black](53)at(2,2){};\node[black](54)at(3,2){};\node[black](55)at(4,2){};
\node[gray] (61)at(0,1){};\node[gray] (62)at(1,1){};\node[gray] (63)at(2,1){};\node[gray] (64)at(3,1){};\node[black](65)at(4,1){};
\node[gray] (71)at(0,0){};\node[black] (72)at(1,0){};\node[black] (73)at(2,0){};\node[black] (74)at(3,0){};\node[black](75)at(4,0){};
\pgfdeclarelayer{bg}    % declare background layer
\pgfsetlayers{bg,main} 
\begin{pgfonlayer}{bg}
\draw (11)--(15);\draw (21)--(25);\draw (31)--(35);\draw (41)--(45);\draw (51)--(55);\draw (61)--(65);\draw(71)--(75);
\draw (11)--(71);\draw (12)--(72);\draw (13)--(73);\draw (14)--(74);\draw (15)--(75);
\end{pgfonlayer}
\node at (-0.3,6){$i_2$};
\node at (-0.3,4){$j_1$};
\node at (4.3,5) {$j_2$};
\node at (-0.3,0) {$i_1$};

\tikzstyle{gline} = [draw = gray!110, line width = 1.7pt]
\tikzstyle{bline} = [draw = black, line width = 1.7pt]
\tikzstyle{tline} = [densely dotted, line cap=round]
\draw[bline](11)--(12)--(13)--(14)--(15)--(25)--(35)--(34)--(33)--(32);
\draw[bline](35)--(45)--(55)--(65)--(75)--(74)--(73)--(72);
\draw[bline](55)--(54)--(53)--(52);
\draw[gline](11)--(21)--(22)--(23)--(24);
\draw[gline](21)--(31)--(41)--(42)--(43)--(44);
\draw[gline](41)--(51)--(61)--(62)--(63)--(64);
\draw[gline](61)--(71);
\end{tikzpicture}
\caption{A hamiltonian coloring with a quartet for $P_5\times P_7$.}\label{fig-odd-grid}
\end{minipage}
\end{figure}
\end{example}

\section{Proofs}\label{sec-proofs}

In this section we prove Theorem \ref{tree}. Figures \ref{fig-basic}, \ref{fig-tree} and \ref{fig-path} illustrate the construction used in the proof. 

\begin{figure}[h]\centering
\begin{minipage}{8cm}
\begin{minipage}{2.5cm}\centering
\begin{tikzpicture}
\tikzstyle{gray}  = [shape = circle, fill = gray!110, draw = gray!110, very thick, minimum size = 2pt,inner sep=2pt]
\tikzstyle{black} = [shape = circle, fill = black, draw = black, very thick , minimum size = 2pt,inner sep=2pt]
\node[gray] (z1) at (0,2) {};
\node[gray] (y1) at (0,1) {};
\node[gray] (x1) at (0,0) {};
\node[black] (c1) at (1,2) {};
\node[black] (b1) at (1,1) {};
\node[black] (a1) at (1,0) {};
\tikzstyle{gline} = [draw = gray!110, line width = 1.5pt]
\tikzstyle{bline} = [draw = black, line width = 1.5pt]
\tikzstyle{tline} = [densely dotted, line cap=round]
%inner lines
\draw[gline](x1)--(z1); \draw[bline](a1)--(c1);\draw[tline](z1)--(b1);\draw[tline](c1)--(y1);
\pgfdeclarelayer{bg}    % declare background layer
\pgfsetlayers{bg,main} 
\begin{pgfonlayer}{bg}
\fill[gray!20] (-0.3,-0.3)--(1.3,-0.3)--(1.3,0.3)--(-0.3,0.3)--(-0.3,-0.3);
\fill[gray!20] (-0.3,1.7)--(1.3,1.7)--(1.3,2.3)--(-0.3,2.3)--(-0.3,2.3);
\end{pgfonlayer}
\node () at (0.6,0.0) {$J$};
\node () at (0.5,1.9) {$I$};
\end{tikzpicture}\\
$C$
\end{minipage}
\begin{minipage}{2.5cm}\centering
\begin{tikzpicture}
\tikzstyle{gray}  = [shape = circle, fill = gray!110, draw = gray!110, very thick, minimum size = 2pt,inner sep=2pt]
\tikzstyle{black} = [shape = circle, fill = black, draw = black, very thick , minimum size = 2pt,inner sep=2pt]
\node[black] (z1) at (0,2) {};
\node[gray] (y1) at (0,1) {};
\node[gray] (x1) at (0,0) {};
\node[gray] (c1) at (1,2) {};
\node[black] (b1) at (1,1) {};
\node[black] (a1) at (1,0) {};
\tikzstyle{gline} = [draw = gray!110, line width = 1.5pt]
\tikzstyle{bline} = [draw = black, line width = 1.5pt]
\tikzstyle{tline} = [densely dotted, line cap=round]
%inner lines
\draw[gline](x1)--(y1)--(c1); \draw[bline](a1)--(b1)--(z1);\draw[tline](z1)--(y1);\draw[tline](c1)--(b1);
\pgfdeclarelayer{bg}    % declare background layer
\pgfsetlayers{bg,main} 
\begin{pgfonlayer}{bg}
\fill[gray!20] (-0.3,-0.3)--(1.3,-0.3)--(1.3,0.3)--(-0.3,0.3)--(-0.3,-0.3);
\fill[gray!20] (-0.3,1.7)--(1.3,1.7)--(1.3,2.3)--(-0.3,2.3)--(-0.3,2.3);
\end{pgfonlayer}
\node () at (0.6,0.0) {$J$};
\node () at (0.5,1.9) {$I$};
\end{tikzpicture}\\
$C\Delta I$
\end{minipage}
\begin{minipage}{2.5cm}\centering
\begin{tikzpicture}
\tikzstyle{gray}  = [shape = circle, fill = gray!110, draw = gray!110, very thick, minimum size = 2pt,inner sep=2pt]
\tikzstyle{black} = [shape = circle, fill = black, draw = black, very thick , minimum size = 2pt,inner sep=2pt]
\node[gray] (z1) at (0,2) {};
\node[gray] (y1) at (0,1) {};
\node[black] (x1) at (0,0) {};
\node[black] (c1) at (1,2) {};
\node[black] (b1) at (1,1) {};
\node[gray] (a1) at (1,0) {};
\tikzstyle{gline} = [draw = gray!110, line width = 1.5pt]
\tikzstyle{bline} = [draw = black, line width = 1.5pt]
\tikzstyle{tline} = [densely dotted, line cap=round]
%inner lines
\draw[gline](y1)--(z1); \draw[bline](b1)--(c1);\draw[tline](z1)--(b1)--(a1);\draw[tline](c1)--(y1)--(x1);
\pgfdeclarelayer{bg}    % declare background layer
\pgfsetlayers{bg,main} 
\begin{pgfonlayer}{bg}
\fill[gray!20] (-0.3,-0.3)--(1.3,-0.3)--(1.3,0.3)--(-0.3,0.3)--(-0.3,-0.3);
\fill[gray!20] (-0.3,1.7)--(1.3,1.7)--(1.3,2.3)--(-0.3,2.3)--(-0.3,2.3);
\end{pgfonlayer}
\node () at (0.6,0.0) {$J$};
\node () at (0.5,1.9) {$I$};
\end{tikzpicture}\\
$C\Delta J$
\end{minipage}
\caption{A hamiltonian coloring $C$ of a graph with a quartet $(I,J)$.}\label{fig-basic}
\end{minipage}$\quad$
\begin{minipage}{5cm}\centering
\begin{tikzpicture}
\tikzstyle{gray}  = [shape = circle, fill = gray!110, draw = gray!110, very thick, minimum size = 2pt,inner sep=2pt]
\tikzstyle{black} = [shape = circle, fill = black, draw = black, very thick , minimum size = 2pt,inner sep=2pt]
\node[black] (a) at (0,1) {};
\node[black] (s) at (1,1) {};
\node[black] (d) at (2,1) {};
\node[black] (f) at (3,1) {};
\node[black] (g) at (4,1) {};
\node[black] (h) at (5,1) {};
\node[black] (q) at (0,0) {};
\node[black] (w) at (1,0) {};
\node[black] (e) at (2,0) {};
\node[black] (r) at (3,0) {};
\node[black] (t) at (4,0) {};
\draw (a)--(s)--(d)--(f)--(g)--(h);
\draw (q)--(w)--(e);
\draw (f)--(r)--(t);
\draw (s)--(w);
\node () at (0,1.3) {$r$};
\node () at (5,1.3) {$l$};
\end{tikzpicture}
\caption{A tree.}\label{fig-tree}
\end{minipage}
\begin{minipage}{14cm}\centering
\begin{tikzpicture}
\node () at (0,3) {}; %this just add some space above
\tikzstyle{gray}  = [shape = circle, fill = gray!110, draw = gray!110, very thick, minimum size = 2pt,inner sep=2pt]
\tikzstyle{black} = [shape = circle, fill = black, draw = black, very thick , minimum size = 2pt,inner sep=2pt]

\node[gray] (z1) at (0,2) {};
\node[gray] (y1) at (0,1) {};
\node[gray] (x1) at (0,0) {};
\node[black] (c1) at (1,2) {};
\node[black] (b1) at (1,1) {};
\node[black] (a1) at (1,0) {};

\node[gray] (z2) at (2.5,2) {};
\node[black] (y2) at (2.5,1) {};
\node[black] (x2) at (2.5,0) {};
\node[black] (c2) at (3.5,2) {};
\node[gray] (b2) at (3.5,1) {};
\node[gray] (a2) at (3.5,0) {};

\node[gray] (z3) at (5,2) {};
\node[gray] (y3) at (5,1) {};
\node[gray] (x3) at (5,0) {};
\node[black] (c3) at (6,2) {};
\node[black] (b3) at (6,1) {};
\node[black] (a3) at (6,0) {};

\node[gray] (z4) at (7.5,2) {};
\node[black] (y4) at (7.5,1) {};
\node[black] (x4) at (7.5,0) {};
\node[black] (c4) at (8.5,2) {};
\node[gray] (b4) at (8.5,1) {};
\node[gray] (a4) at (8.5,0) {};

\node[gray] (z5) at (10,2) {};
\node[gray] (y5) at (10,1) {};
\node[gray] (x5) at (10,0) {};
\node[black] (c5) at (11,2) {};
\node[black] (b5) at (11,1) {};
\node[black] (a5) at (11,0) {};

\node[gray] (z6) at (12.5,2) {};
\node[black] (y6) at (12.5,1) {};
\node[black] (x6) at (12.5,0) {};
\node[black] (c6) at (13.5,2) {};
\node[gray] (b6) at (13.5,1) {};
\node[gray] (a6) at (13.5,0) {};

\node[gray]  (lz1) at (0,-2) {};
\node[black] (ly1) at (0,-3) {};
\node[black] (lx1) at (0,-4) {};
\node[black] (lc1) at (1,-2) {};
\node[gray]  (lb1) at (1,-3) {};
\node[gray]  (la1) at (1,-4) {};

\node[gray]  (lz2) at (2.5,-2) {};
\node[gray]  (ly2) at (2.5,-3) {};
\node[gray]  (lx2) at (2.5,-4) {};
\node[black] (lc2) at (3.5,-2) {};
\node[black] (lb2) at (3.5,-3) {};
\node[black] (la2) at (3.5,-4) {};

\node[gray]  (lz3) at (5,-2) {};
\node[black] (ly3) at (5,-3) {};
\node[black] (lx3) at (5,-4) {};
\node[black] (lc3) at (6,-2) {};
\node[gray]  (lb3) at (6,-3) {};
\node[gray]  (la3) at (6,-4) {};

\node[gray]  (lz4) at (7.5,-2) {};
\node[gray]  (ly4) at (7.5,-3) {};
\node[gray]  (lx4) at (7.5,-4) {};
\node[black] (lc4) at (8.5,-2) {};
\node[black] (lb4) at (8.5,-3) {};
\node[black] (la4) at (8.5,-4) {};

\node[gray]  (lz5) at (10,-2) {};
\node[black] (ly5) at (10,-3) {};
\node[black] (lx5) at (10,-4) {};
\node[black] (lc5) at (11,-2) {};
\node[gray]  (lb5) at (11,-3) {};
\node[gray]  (la5) at (11,-4) {};

\tikzstyle{gline} = [draw = gray!110, line width = 1.5pt]
\tikzstyle{bline} = [draw = black, line width = 1.5pt]
\tikzstyle{tline} = [densely dotted, line cap=round]

%inner lines
\draw[gline](x1)--(z1); \draw[bline](a1)--(c1);\draw[tline](z1)--(b1);\draw[tline](c1)--(y1);
\draw[bline](x2)--(y2)--(c2); \draw[gline](a2)--(b2)--(z2);\draw[tline](y2)--(z2);\draw[tline](b2)--(c2);

\draw[gline](x3)--(z3); \draw[bline](a3)--(c3);\draw[tline](z3)--(b3);\draw[tline](c3)--(y3);
\draw[bline](x4)--(y4)--(c4); \draw[gline](a4)--(b4)--(z4);\draw[tline](y4)--(z4);\draw[tline](b4)--(c4);

\draw[gline](x5)--(z5); \draw[bline](a5)--(c5);\draw[tline](z5)--(b5);\draw[tline](c5)--(y5);
\draw[bline](x6)--(y6)--(c6); \draw[gline](a6)--(b6)--(z6);\draw[tline](y6)--(z6);\draw[tline](b6)--(c6);

%lower inner lines
\draw[bline](lx1)--(ly1)--(lc1); \draw[gline](la1)--(lb1)--(lz1);\draw[tline](ly1)--(lz1);\draw[tline](lb1)--(lc1);

\draw[gline](lx2)--(lz2); \draw[bline](la2)--(lc2);\draw[tline](lz2)--(lb2);\draw[tline](lc2)--(ly2);

\draw[bline](lx3)--(ly3)--(lc3); \draw[gline](la3)--(lb3)--(lz3);\draw[tline](ly3)--(lz3);\draw[tline](lb3)--(lc3);

\draw[gline](lx4)--(lz4); \draw[bline](la4)--(lc4);\draw[tline](lz4)--(lb4);\draw[tline](lc4)--(ly4);

\draw[bline](lx5)--(ly5)--(lc5); \draw[gline](la5)--(lb5)--(lz5);\draw[tline](ly5)--(lz5);\draw[tline](lb5)--(lc5);

%upper arcs
\draw[gline](z2) to[bend left] (z3);
\draw[gline](z3) to[bend left] (z4);
\draw[gline](z4) to[bend left] (z5);
\draw[gline](z5) to[bend left] (z6);
\draw[bline](c2) to[bend left] (c3);
\draw[bline](c3) to[bend left] (c4);
\draw[bline](c4) to[bend left] (c5);
\draw[bline](c5) to[bend left] (c6);
\draw[gline](z1) to[bend left] (z2);
\draw[bline](c1) to[bend left] (c2);
%lower arcs
\draw[tline](x1) to[bend right] (x2);
\draw[tline](x2) to[bend right] (x3);
\draw[tline](a1) to[bend right] (a2);
\draw[tline](a2) to[bend right] (a3);
\draw[tline](a3) to[bend right] (a4);
\draw[tline](a4) to[bend right] (a5);
\draw[tline](a5) to[bend right] (a6);
\draw[tline](x3) to[bend right] (x4);
\draw[tline](x4) to[bend right] (x5);
\draw[tline](x5) to[bend right] (x6);
%middle arcs
\draw[tline](b1) to[bend right] (b2);
\draw[tline](b2) to[bend right] (b3);
\draw[tline](b3) to[bend right] (b4);
\draw[tline](b4) to[bend right] (b5);
\draw[tline](b5) to[bend right] (b6);
\draw[tline](y1) to[bend right] (y2);
\draw[tline](y2) to[bend right] (y3);
\draw[tline](y3) to[bend right] (y4);
\draw[tline](y4) to[bend right] (y5);
\draw[tline](y5) to[bend right] (y6);
%lower upper arcs
%upper arcs
\draw[gline](lz2) to[bend left] (lz3);
\draw[gline](lz4) to[bend left] (lz5);
\draw[bline](lc2) to[bend left] (lc3);
\draw[bline](lc4) to[bend left] (lc5);
\draw[gline](lz1) to[bend left] (lz2);
\draw[bline](lc1) to[bend left] (lc2);
%lower lower arcs
\draw[tline](lx1) to[bend right] (lx2);
\draw[tline](lx2) to[bend right] (lx3);
\draw[tline](la1) to[bend right] (la2);
\draw[tline](la2) to[bend right] (la3);
\draw[tline](la4) to[bend right] (la5);
\draw[tline](lx4) to[bend right] (lx5);
%middle arcs
\draw[tline](lb1) to[bend right] (lb2);
\draw[tline](lb2) to[bend right] (lb3);
\draw[tline](lb4) to[bend right] (lb5);
\draw[tline](ly1) to[bend right] (ly2);
\draw[tline](ly2) to[bend right] (ly3);
\draw[tline](ly4) to[bend right] (ly5);
%vertical lines
\draw[gline](z2) to[bend right] (lz2);
\draw[tline](y2) to[bend right] (ly2);
\draw[tline](x2) to[bend right] (lx2);
\draw[bline](c2) to[bend left] (lc2);
\draw[tline](b2) to[bend left] (lb2);
\draw[tline](a2) to[bend left] (la2);

\draw[gline](z4) to[bend right] (lz4);
\draw[tline](y4) to[bend right] (ly4);
\draw[tline](x4) to[bend right] (lx4);
\draw[bline](c4) to[bend left] (lc4);
\draw[tline](b4) to[bend left] (lb4);
\draw[tline](a4) to[bend left] (la4);

\pgfdeclarelayer{bg}    % declare background layer
\pgfsetlayers{bg,main} 
\begin{pgfonlayer}{bg}
\fill[gray!30] (-0.3,-0.3)--(1.3,-0.3)--(1.3,0.3)--(-0.3,0.3)--(-0.3,-0.3);
\fill[gray!30] (12.3,1.7)--(13.8,1.7)--(13.8,2.3)--(12.3,2.3)--(12.3,2.3);
\end{pgfonlayer}

\node () at (0.6,0.0) {$J_r$};
\node () at (13,1.9) {$I_l$};

\end{tikzpicture}
\caption{A construction for the hamiltonian coloring $D$ of $G\times T$ as in the proof of Theorem \ref{tree} with quartet $(J_r,I_l)$, where $G$ is the graph in Figure \ref{fig-basic} and $T$ is the tree in Figure \ref{fig-tree}.}
\label{fig-path}
\end{minipage}
\end{figure}

%Let $G$ be a graph and suppose that $V(G):=X_1\dot{\cup}\cdots\dot{\cup}X_n$. Let $K$ be a set of two colors. For $i=1,\dots,n$, let $C_i$ is a coloring of $X_i$ with colors in $K$. We denote by $C_1\cup\cdots C_n$ the coloring of $V(G)$ such that for each $v\in V(G)$, the color of $v$, is $C_i(v)$, where $i$ is the index satisfying $v\in X_i$. 
\bigskip
\newcommand{\ired}{i_{\rm red}}
\newcommand{\iblue}{i_{\rm blue}}
\newcommand{\red}{{\rm red}}
\newcommand{\blue}{{\rm blue}}

Consider a graph $G$ and a coloring $C:G\rightarrow\{\red,\blue\}$.  We say that an edge of $G$ is a \defin{red} (resp. \defin{blue}) edge of $C$ if its both endvertices have color red (resp. blue) in $C$. An edge of $G$ is \defin{monochromatic} if it is red or blue. \\

\begin{proofof}\emph{Proof of Theorem \ref{tree}: }
We may assume that $|T|\ge2$. For each $x\in V(T)$, we denote by $G_x$ the graph $G\times T[\{x\}]$ and, for $A\cont V(G)$, $A_x:=A\times \{x\}$.

Let $C$ be a hamiltonian coloring of $G$ with a quartet $(I,J)$ and colors in $\{\red,\blue\}$. We write $I=\{i_\red,i_\blue\}$ and $J=\{j_\red,j_\blue\}$, where the vertices are colored by $C$ with the color in their indices.

For each $x\in V(T)$, let $C_x$ be the coloring of $G_x$ such that $C_x(v,x):=C(v)$ for all $v\in V(G)$. Note that $C_x$ is a hamiltonian coloring of $G_x$ with quartet $(I_x,J_x)$.

Consider a leaf $r$ as root for $T$ and, for each $x\in V(T)$, let $d(x)$ denote the usual distance from $r$ to $x$ in $T$. Now, for each $x\in V(T)$, define $D_x:=C_x$ if $d(x)$ is even and $D_x:=C_x\Delta(V(G_x)-I_x)$ if $d(x)$ is odd. Let $D$ be the $2$-coloring of $G\times T$ such that, for each $x\in V(T)$, the restriction of $D$ to $G_x$ is equal to $D_x$. Note that for all $x\in V(T)$, $D(\ired,x)={\rm red}$ and $D(\iblue,x)={\rm blue}$. Let $l$ be a leaf of $T$ other than $r$. We shall prove that $D$ is a hamiltonian coloring of $G\times T$ with $(J_r,I_l)$ as quartet. Next we check:\\

\noindent{\bf Claim 1.} {\it For each $\kappa\in\{\red,\blue\}$ and each $xy\in E(T)$, $(i_\kappa,x)(i_\kappa,y)$ is the unique edge of color $\kappa$ linking a vertex of $G_x$ to a vertex of $G_y$.}\\

Indeed, by the definition of $G\times T$, such an edge must be of the form $(v,x)(v,y)$ with $v\in V(G)$. We may assume without loosing generality that $d(x)$ is even and $d(y) $ is odd. So $D_x=C_x$ and $D_y=C_y\Delta (V(G_y)-I_y)$. Thus the edge $(v,x)(v,y)$ is monochromatic if and only if $v\in I$. Moreover, this edge has color $\kappa$ if and only if $C(v)=\kappa$ and this implies Claim 1.\\

Denote $H:=G\times T$ and let $H_\red$ be the subgraph of $H$ induced by the red vertices of $D$. Since $T$ is a tree, Claim 1 implies:\\

\noindent{\bf Claim 2.} {\it Let $xy\in E(G)$, then $(\ired,x)(\ired,y)$ is a bridge of $H_\red$.}\\

By Claim 2, if $Z$ is a cycle of $R$, then $Z$ is a cycle of $G_x$ for some $x\in V(T)$. But, as $C_x$ is a hamiltonian coloring of $G_x$, it follows by (Q2) that $D_i$ is also a hamiltonian coloring of $G_x$. Thus $H_\red\cap G_x$ is a tree and $H_\red$ has no cycles. Note $G[\{\ired\}]\times T$ is a connected subgraph of $H_\red$ because it is isomorphic to $T$. Each vertex of $H_\red$ is in a subgraph of the form $H_\red \cap G_x$, which is a tree with a vertex in common with $G[\{\ired\}]\times T$. Thus $H_\red$ is connected and, therefore, $H_\red$ is a tree. Analogously the same holds for the subgraph of $H$ induced by the blue vertices of $D$. This proves that $D$ is a hamiltonian coloring of $H$.

It is left to prove that $(J_r,I_l)$ is a quartet for $D$. For this we have to prove (Q1)-(Q3). By construction, both $J_l$ and $I_r$ have vertices of different colors. So (Q1) holds.

%For each coloring $K:V(G\times T)\rightarrow\{\red,\blue\}$, we denote by $R(K):=\{(v,x)\in V(G\times T):K(v,x)=\red\}$. 

Consider the coloring $D':=D\Delta J_r$. Denote by $R'$ the set of red vertices of $D'$ and let $H'_\red:=H[R']$. Note that $H_1:=H'_\red-(V(G_r)\cap V(H'_\red))=H_\red-(V(G_r)\cap V(H_\red))$. Since $H_\red$ is a tree, it follows from Claim 2 that $H_1$ is a tree. Now consider the restriction $D'_x$ of $D'$ to $G_x$. As $D_i=C_i$ is a hamiltonian coloring of $G_r$ with $(I_r,J_r)$ as quartet, it follows from (Q3) that the red vertices of $D'_r=C_r\Delta J_r$ induce in $G_x$ a forest with exactly two components $T_I$ and $T_J$ such that $V(T_I)\cap I_r=\{(i_\red,r)\}$ and $V(T_J)\cap J_r=\{(j_\blue,r)\}$. Let $s$ be the neighbor of $r$ in $T$. It follows from Claim 1 that $(\ired,r)(\ired,s)$ and $(j_\blue,r)(j_\blue,s)$ are the unique edges that link $G_x\cap H'_\red$ to $H_1$. This implies that $H'_\red$ is obtained from the disjoint trees $H_1$, $T_J$ and $T_J$ by adding $(\ired,r)(\ired,s)$ linking $T_i$ to $H_1$ and adding $(j_\blue,r)(j_\blue,s)$ linking $T_j$ to $H_1$ and, therefore, $H'_\red$ is a tree. Analogously, the blue vertices of $D'$ also induce a tree in $H$ and this proves (Q2).

Consider now $D'':=D\Delta I_l$. Let $R''$ be the set of red vertices of $D''$ and $H''_\red:=H[R'']$. Similarly to $H_1$, $H_2:=H''_\red-(V(G_l)\cap V(H''_\red))=H_\red-(V(G_l)\cap V(H_\red))$ is a tree. By Claim 1 there is no monochromatic edge of $D''$ from $H_2$ to $G_l$. The restriction $D''_l$ of $D''$ to $G_l$ is equal to $C_l\Delta I_l$ if $d(l)$ is even and equal to $C_l\Delta V(G_l)$ if $d(l)$ is odd. Thus $D''_l$ is a hamiltonian coloring of $G_l$ and $H''_\red\cap G_l$ is a tree. Therefore, $H''_\red$ is a forest with $H_2$ and $H''_\red\cap G_l$ as connected components. Note that $H_2$ has a vertex of $J_r$ and $H''_\red\cap G_l$ has a vertex of $I_l$. Since similar arguments are valid for the blue vertices, this implies that (Q3) holds and the theorem is valid.
\end{proofof}

\section*{Acknowledgments}

 The Author thanks to Haindong Wu for the helpful discussions on the subject.

%\providecommand{\bysame}{\leavevmode\hbox to3em{\hrulefill}\thinspace}
%\providecommand{\MR}{\relax\ifhmode\unskip\space\fi MR }
%% \MRhref is called by the amsart/book/proc definition of \MR.
%\providecommand{\MRhref}[2]{%
%  \href{http://www.ams.org/mathscinet-getitem?mr=#1}{#2}
%}
%\providecommand{\href}[2]{#2}

\end{document}